\documentclass[12pt]{article}
\usepackage[cp1251]{inputenc}
\usepackage[russian]{babel}
\usepackage{amsmath,amsthm,amssymb}
\usepackage{latexsym}
\usepackage{amsfonts}
\newcommand{\en}{\enspace}

\begin{document}

$$
$$

\begin{center}
\bf COMPLEX POTENTIALS: BOUND STATES,\\
QUANTUM DYNAMICS AND WAVE OPERATORS
\end{center}

\medskip

\begin{center}
\bf S. A. Stepin
\end{center}

\bigskip
\bigskip

\textbf{Abstract.}
Schr\"odinger operator on half-line with complex potential and the corresponding evolution are studied within perturbation theoretic approach. The total number of eigenvalues and spectral singularities is effectively evaluated. Wave operators are constructed and a criterion is established for the similarity of perturbed and free propagators.

\medskip

\textbf{Keywords:} Schr\"odinger operator, propagator, Jost function, wave operators

\medskip

{\bfseries 2010 Mathematics Subject Classification:} 34L15, 34L25

\bigskip
\bigskip

{\bf 1. Introduction}

\medskip

Given bounded complex valued potential $\,V(x)\,$ consider in $\,{\mathcal H}={\rm L}_2(0,\infty)\,$ an operator
$$
L_V\,=\,L_0\,+\,V\,=\,-d^2/dx^2\,+\,V(x)
$$
generated by Dirichlet boundary condition at zero. Such an operator proves to be quite a simple and rather capacious model which displays a number of effects typical for perturbation theory in nonselfadjoint setting (see [1] and [2]). Besides that, Schr\"odinger operator with complex potential is known (see [3]) to appear in the study of open quantum mechanical systems with energy dissipation.

Certain difficulties in the study of Schr\"odinger operators with complex potential are due to the lack of an a priori information about the behavior of the resolvent $\,R_V(\lambda)=(L_V-\lambda I)^{-1}\,$ near the points of the spectrum $\,\sigma(L_V)\,$ as well as complicated structure of the spectrum itself. In this context an essential role is played (see [4]) by the spectral singularities, i.e. poles of the analytic continuation of the resolvent integral kernel which are embedded into continuous spectrum. It turns out that the crucial information about spectral properties of operator $\,L_V\,$ can be extracted from the so-called Jost function which coincides under appropriate assumptions (see [5]) with Fredholm determinant
$$
e(k)\en=\en\det\big(\,I\,+\,V\,|V|^{-1/2}(L_0-k^2I)^{-1}|V|^{1/2}\big)\,.
$$

For one-dimensional Schr\"odinger operators the so-called transformation operators prove to be an effective tool in the study of spectral similarity as well as for the solving of direct and inverse scattering problems. In turn scattering theory itself provides an adequate construction of transformation operators intertwining perturbed operator with unperturbed one. This approach was elaborated in [6] and [7] for Schr\"odinger operator with complex potential; to this end wave operators are constructed by means of comparison of the corresponding perturbed and unperturbed propagators $\,e^{itL_V}\,$ and $\,e^{itL_0}.\,$ Note that the problem concerning similarity of the part of operator $\,L_V\,$ associated with its continuous spectrum and unperturbed operator $\,L_0\,$ corresponding to $\,V(x)\equiv 0\,$ is closely related to generalized eigenfunction expansion problem (cf. [8]).

The present paper is organized as follows. In Section 1 some auxiliary estimates for certain solutions to stationary Schr\"odinger equation are derived. These estimates will be used both in the proofs of Theorem 1 and Theorem 2 below.

In Section 2 we study discrete spectrum of operator  $\,L_V.\,$ Theorem 1 gives an effective upper bound for the total number of both its eigenvalues and spectral singularities under certain assumptions imposed on the potential $\,V(x)\,$ which allow the resolvent integral kernel to possess an analytic continuation through the continuous spectrum cut. Such estimates in nonselfadjoint case have been previously obtained for Boltzmann transport operator in [9] (see also [10]).

Analytic properties of the resolvent $\,R_V(\lambda)=(L_V-\lambda I)^{-1}\,$ are investigated in Section 3.
Within the present context the notion of relative smoothness which is due to Kato will be useful. Operator $\,A\,$ is said to be smooth with respect to $\,L_V\,$ if for an arbitrary $\,\varphi\in\mathcal H\,$ vector-valued function $\,AR_V(\lambda)\varphi\,$ belongs to Hardy classes $\,{\rm H}_2^{\pm}\,$ in upper and lower half-planes $\,\mathbb C_{\pm}.\,$ For potentials $\,V(x)\,$ such that
$$
\int_0^{\infty}x|V(x)|\,dx\,<\,\infty \eqno (\ast)
$$
operator $\,A=\sqrt{V}\,$ proves to be smooth with respect to $\,L_V\,$ provided that spectrum $\,\sigma(L_V)\,$ is purely continuous without spectral singularities embedded.

In Section 4 we prove Theorem 2 which is a criterion of similarity $\,L_V\sim L_0\,$ within the class of potentials possessing the first momentum $(\ast).$ It was established in [6] that under the condition
$$
\int_0^{\infty}x|V(x)|\,dx\,<\,1\eqno(\ast\ast)
$$
operator $\,L_V=L_0+V\,$ is similar to $\,L_0\,$ and moreover this similarity is implemented by explicit construction of the
corresponding wave operators. Sufficient condition $(\ast\ast)$ is sharp in the sense that numerical upper bound is the best possible. Otherwise a counterexample shows that an obstacle to the similarity is emission of an eigenvalue of $\,L_V.\,$
Theorem 2 thus extends and supplements the class of potentials restricted by the Kato condition $(\ast\ast).$

\bigskip

{\bf 2. Estimates for solutions to Schr\"odinger equation}

\medskip

Provided that potential $\,V(x)\,$ is integrable on half-line $\,\mathbb R_+=[0,\infty)\,$ the corresponding  Schr\"odinger equation
\begin{equation}
-y''\,+\,V(x)y\,=\,k^2y
\label{form1}
\end{equation}
for $\,k\in\mathbb C_+\,$ is known to have (see [1]) the so-called Jost solution $\,e(x,k)\,$ specified by its asymptotics $\,e(x,k)\sim e^{ikx}\,$ at $\infty.$

\bigskip

{\bf Lemma 1.} {\it Suppose that for a certain $\,a\geqslant 0\,$ potential $\,V(x)\,$ satisfies condition
\begin{equation}
\int_0^{\infty}e^{ax}|V(x)|\,dx\,<\,\infty\,.
\end{equation}
Then for arbitrary $\,k\in\mathbb C,\, {\rm Im}\,k>-a/2,\,$ and $\,\alpha\in[0,1]\,$ Jost solution to equation {\rm (1)} admits for $\,x\geqslant 0\,$ the estimate
$$
|\,e(x,k)e^{-ikx}-1|\,\leqslant\,\exp\left(\frac1{(2|k|)^{1-\alpha}}\int_x^{\infty}\!\!\xi^{\alpha}\max\{1,e^{-2{\rm Im}\,k\xi}\}^{\alpha}\big(1+e^{-2{\rm Im}\,k\xi}\big)^{1-\alpha}|V(\xi)|\,d\xi\right)-\,1.
$$
}

For $\,k\in\mathbb C\,$ such that $\,{\rm Im}\,k>-a/2\,$ set $\,\varepsilon^{(0)}(x,k)=1\,$ and
$$
\varepsilon^{(n+1)}(x,k)\,=\,\int_x^{\infty}\frac{e^{2ik(\xi-x)}-1}{2ik}\,V(\xi)\,\varepsilon^{(n)}(\xi,k)\,d\xi\,.
$$
Making use of the inequality $\,|\,e^{iz}-1|\leqslant|z|\max\{1,e^{-{\rm Im}\,z}\},\,$ one can inductively verify the estimate
$$
|\,\varepsilon^{(n)}(x,k)|\,\leqslant\,\frac1{n!}\left(\frac1{(2|k|)^{1-\alpha}}\int_x^{\infty}\!\!\xi^{\alpha}\max\{1,e^{-2{\rm Im}\,k\xi}\}^{\alpha}\big(1+e^{-2{\rm Im}\,k\xi}\big)^{1-\alpha}|V(\xi)|\,d\xi\right)^n
$$
valid for arbitrary $\,\alpha\in [0,1].\,$
Due to this fact under the condition (2) series $\,e(x,k)=e^{ikx}\sum\limits_{n=0}^{\infty}\varepsilon^{(n)}(x,k)\,$ converges uniformly in $\,k\,$ and $\,x,\,$ hence it represents a solution to integral equation
$$
e(x,k)\,=\,e^{ikx}\,-\,\int_x^{\infty}\frac{\sin k(x-\xi)}{k}\,V(\xi)\,e(\xi,k)\,d\xi
$$
being exactly Jost solution to (1) appropriately evaluated.

\bigskip

{\bf Corollary 1.} {\it In particular when $\,k\in\mathbb C_+\,$ one gets the following estimate
\begin{equation}
|\,e(x,k)|\en\leqslant\en\exp\left(\int_x^{\infty}\xi|V(\xi)|\,d\xi\right)\,e^{-{\rm Im}\,k\,x}\,.
\end{equation}
}

Denote by $\,s(x,k)\,$ a solution to equation (1) determined by initial data $\,s(0,k)=0\,$ and $\,s'_x(0,k)=1.\,$
It proves to be a solution of integral equation
\begin{equation*}
s(x,k)\,\,=\,\,\frac{\sin kx}k\,+\,\int_0^x\frac{\sin
k(x-\xi)}{k}\,V(\xi)\,s(\xi,k)\,d\xi\,
\end{equation*}
and hence the following inequality
$$
|s(x,k)|\,\leqslant\,xe^{{\rm Im}\,k\,x}\bigg( 1\,+\,\int_0^x
e^{-{\rm Im}\,k\,\xi}|V(\xi)||s(\xi,k)|\,d\xi\bigg)
$$
holds true since $\,|\sin
kx|\leqslant |k|\,xe^{{\rm Im}\,k\, x}.\,$ Applying Gronwall's lemma one immediately gets

\medskip

{\bf Lemma 2.} {\it Provided that $\,\displaystyle{\int_0^{\infty}\!x|V(x)|\,dx<\infty}\,$ the estimate
\begin{equation}
|s(x,k)|\,\leqslant\, x\,\exp\bigg(\int_0^x\xi|V(\xi)|\,d\xi\bigg)e^{{\rm Im}\,k\,x}
\end{equation}
is valid for arbitrary $\,k\in\mathbb C_+\,$ and all $\,x\geqslant 0.$
}

\bigskip

{\bf 3. Estimate for the number of bound states}

\medskip

The spectrum of operator $\,L_V=L_0+V\,$ with complex-valued potential $\,V(x)\,$ integrable on $\,\mathbb R_+\,$
consists of continuous and discrete components
$$
\sigma_c(L_V)=\,\mathbb R_+,\quad \sigma_d(L_V)\,=\,\{k^2:e(k)=0,k\in\mathbb C_+\}\,,
$$
where $\,e(k):=e(0,k)\,$ is the so-called Jost function. By virtue of Lemma 1 Jost function is analytic in $\,\mathbb C_+,\,$ it admits continuation to $\,\mathbb R\setminus\{0\}\,$ and satisfies the estimate
\begin{equation}
|\,e(k)-1|\leqslant\exp\left(\frac1{(2|k|)^{1-\alpha}}\int_0^{\infty}\!\!\xi^{\alpha}\big(1+e^{-2{\rm Im}\,k\xi}\big)|V(\xi)|\,d\xi\right)-1.
\end{equation}
Thus the set $\,\sigma_d(L_V)\,$ of eigenvalues of operator $\,L_V\,$ is bounded, at most countable and its accumulation points (if any) belong to the half-line $\,\mathbb R_+.\,$ Note that operator $\,L_V\,$ has no positive and thus embedded eigenvalues.

\medskip

{\bf Definition.} {\it Real zeroes of Jost function $\,e(k)\,$ correspond to the distinguished points $\,\lambda=k^2\,$ of continuous spectrum $\,\sigma_c(L_V)\,$ called spectral singularities.}

\medskip

Provided that $\,\displaystyle{\int_0^{\infty}\!\!x|V(x)|\,dx<\infty}\,$ Jost function is known to be continuous up to the real line everywhere including zero. Therefore one has the following (cf. [11])

\medskip

{\bf Statement 1.} {\it Operator $\,L_V\,$ without spectral singularities has finite discrete spectrum.}

\medskip

Denote by $\,R\,$ the minimal radius of discs containing all the zeroes of Jost function in the closed upper half-plane. Estimate (5) implies that
$$
R\en\leqslant\en \widetilde{R}(\alpha)\,:=\,\left(\frac1{\ln 2}\int_0^{\infty}(2x)^{\alpha}|V(x)|\,dx\right)^{1/(1-\alpha)}\!,\quad \alpha\in[0,1)\,.
$$
An upper bound for the number (i.e. total multiplicity) $\,N(V)\,$ of eigenvalues and spectral singularities of operator $\,L_V=L_0+V\,$ is given by

\bigskip

{\bf Theorem 1.}\,{\it Suppose that for a certain $\,a>0\,$ integral {\rm (2)} converges. Then given arbitrary $\,\alpha,\beta\in[0,1)\,$ and $\,A>\max\{\widetilde{R}(\alpha),R^2/a-a/4\}\,$ one has the inequality
\begin{equation*}
N(V)\en\leqslant\en\left(\ln\frac{A+a/2}{\sqrt{A^2+R^2}}\right)^{-1}\!\!\bigg\{\frac1{a^{1-\beta}}
\int_0^{\infty}\!x^{\beta}(1+e^{ax})|V(x)|\,dx\,\,-\,\,\ln\Big(2-2^{(\widetilde{R}/A)^{1-\alpha}}\Big)\bigg\}.
\end{equation*}
}

\smallskip

{\bf Proof.} Total multiplicity $\,N(V)\,$ of eigenvalues and spectral singularities of operator $\,L_V\,$ coincides with the number of zeroes of Jost function $\,e(k)\,$ in the closed upper half-plane $\,\overline{\mathbb C}_+.\,$ Moreover all of them are located in the domain
$$
\{k\in\overline{\mathbb C}_+:\,|k|\leqslant R\}\,\subset\,\{|k-iA|\leqslant\sqrt{A^2+R^2}\}.
$$
To estimate the value $\,N(V)\,$ which does not exceed the number of zeroes of function $\,\varphi(z):=e(z+iA)\,$ in the disc $\,|z|\leqslant\sqrt{A^2+R^2}\,$ we apply Nevanlinna-Jensen formula. Note that $\,\sqrt{A^2+R^2}<A+a/2\,$ because $\,A>R^2/a-a/4.\,$ Choosing arbitrary $\,\rho\in(\sqrt{A^2+R^2},A+a/2)\,$ one has
\begin{equation*}
N(V)\,\ln\frac{\rho}{\sqrt{A^2+R^2}}\,\,\,\leqslant\!\sum_{|z_k|^2\leqslant A^2+R^2}\frac{\rho}{|z_k|}\,\,=\,\,
\frac1{2\pi}\int_0^{2\pi}\ln|\varphi(\rho e^{i\theta})|\,d\theta\,-\,\ln|\varphi(0)|.
\end{equation*}
In order to evaluate absolute value of $\,\varphi\,$ on the circle $\,|z|=\rho\,$ use inequality (5) with $\,\alpha\,$ replaced by $\,\beta\,:$
\begin{multline*}
|\varphi(\rho e^{i\theta})|\,=\,|e(iA+\rho e^{i\theta})|\en\leqslant\\
\leqslant\en
\exp\left(\frac1{(2|iA+\rho e^{i\theta}|)^{1-\beta}}\int_0^{\infty}
\!\!\xi^{\beta}\big(1+e^{-2(A+\rho\sin\theta)\xi}\big)|V(\xi)|\,d\xi\right)\\
\leqslant\en\exp\left(\frac1{(2(\rho-A))^{1-\beta}}\int_0^{\infty}
\!\!\xi^{\beta}\big(1+e^{2(\rho-A)\xi}\big)|V(\xi)|\,d\xi\right)\,,
\end{multline*}
Also by virtue of (5) for $\,\alpha\in[0,1)\,$ we obtain
$$
|\varphi(0)|\,\geqslant\,1-|e(iA)-1|\,\geqslant\,2\,-\,
\exp\left[\frac1{A^{1-\alpha}}\int_0^{\infty}\!\!(2x)^{\alpha}
|V(x)|\,dx\right]=\,2-2^{(\widetilde{R}/A)^{1-\alpha}},
$$
where the right-hand side is positive since $\,A>\widetilde{R}(\alpha).\,$ Thus for an arbitrary $\,\rho\in(\sqrt{A^2+R^2},A+a/2)\,$ the inequality
\begin{equation*}
N(V)\,\ln\frac{\rho}{\sqrt{A^2+R^2}}\,\,\leqslant\,\,\frac1{(2(\rho-A))^{1-\beta}}\int_0^{\infty}
\!\!\xi^{\beta}\big(1+e^{2(\rho-A)\xi}\big)|V(\xi)|\,d\xi\,\,-\,\,\ln\Big(2-2^{(\widetilde{R}/A)^{1-\alpha}}\Big)
\end{equation*}
is valid. Passing here to the limit as $\,\rho\to A+a/2\,$ we complete the proof.

\bigskip

{\bf Corollary 2.} {\it
Provided that $\,\,\displaystyle{a\,\geqslant\,b\,=\,\frac1{\ln 2}\int_0^{\infty}\!|V(x)|\,dx}\,\,$
one has the following estimate
\begin{equation*}
N(V)\en\leqslant\en 10\,\bigg\{1\,+\,\,\frac2{b}\int_0^{\infty}\!e^{bx}|V(x)|\,dx\bigg\}.
\end{equation*}
}

\medskip

{\bf 4. Kato smoothness property}

\medskip

Resolvent $\,R_V(\lambda)=(L_V-\lambda I)^{-1}\,$ for
$\,\lambda=k^2,\, k\in\mathbb C_+,\,$ is known (see e.g. [1]) to be an integral operator with the kernel
$$
R_V(x,\xi,\lambda)\,\,=\,\,s(\min\{x,\xi\},k)\,e(\max\{x,\xi\},k)\big/e(k)\,.
$$
Denote by $\,A\,$ and $\,B\,$ operators of multiplication by functions $\,a(x)\,$ and $\,b(x)\,$ such that
$$
\langle\,a\rangle\,:=\,\bigg(\int_0^{\infty}x|a(x)|^2\,dx\bigg)^{1/2}\!<\,\infty\,,\quad
\langle\,b\rangle\,:=\,\bigg(\int_0^{\infty}x|b(x)|^2\,dx\bigg)^{1/2}\!<\,\infty\,.
$$

\medskip

{\bf Lemma 3.}
{\it Under the condition $(\ast)$ operator function
$\,e(k)AR_V(k^2)B:{\mathcal H}\to{\mathcal H}\,$  is analytic in $\,\mathbb C_+\,$ and for all $\,k\in\mathbb C_+\,$ the inequality
$$
\|\,e(k)AR_V(k^2)B\|\,\leqslant\,C\langle\,a\rangle\langle\,b\rangle
$$
holds true with the constant $\,C\,=\,\exp\langle\sqrt{V}\rangle^2.$
}

\bigskip

By virtue of (3) and (4) the integral kernel of the resolvent $\,R_V(\lambda)\,$ satisfies the estimate
$$
|R_V(x,\xi,\lambda)|\,\leqslant\,\frac{C}{|e(k)|}\,\min\{x,\xi\}
$$
and hence one has
\begin{multline*}
\|e(k)AR_V(k^2)Bf\|^2\,\leqslant\,C^2\int_0^{\infty}\!|a(x)|^2
\bigg|\int_0^{\infty}\!\min\{x,\xi\}b(\xi)f(\xi)\,d\xi\,\bigg|^2dx\,\,\leqslant\\
\leqslant\,\,C^2\int_0^{\infty}\!|a(x)|^2
\bigg(\int_0^{\infty}\!(\,\min\{x,\xi\})^2|b(\xi)|^2d\xi\bigg)\bigg(\int_0^{\infty}|f(\xi)|^2d\xi\bigg)dx\,\,\leqslant
\,\, C^2\langle\,a\rangle^2\langle\,b\rangle^2\|f\,\|^2.
\end{multline*}
Resolvent $\,R_V(k^2)\,$ is meromorphic in $\,\mathbb C_+\,$ so that its poles are located at the zeroes of $\,e(k).\,$ In due turn operator function $\,e(k)AR_V(k^2)B\,$ is holomorphic and (according to the above estimate) bounded in the vicinity of each resolvent pole and therefore all these singularities are removable.

\bigskip

{\bf Statement 2.} {\it Suppose that condition $(\ast)$ is fulfilled. If $\,e(k)\ne 0\,$
for $\,k\in\overline{\mathbb C}_+\,$ then operator $\,A=\sqrt{V}\,$ is relatively smooth in the sense of Kato with respect to $\,L_V,\,$ i.e.
\begin{equation}
\sup_{\|\varphi\|=1,\,\varepsilon\geqslant0}\int_{\mathbb
R}\|AR_V(\tau\pm
i\varepsilon)\varphi\|^2\,d\tau\,<\,\infty\,.
\end{equation}
}

\medskip

Indeed, according to resolvent identity one has
\begin{equation}
AR_V(\lambda)\,=\,\big(I-AR_V(\lambda)B\big)AR_0(\lambda)\,
\end{equation}
where $\,B=\sqrt{|V|}.\,$ By Lemma 3 operator $\,BR_0(k^2)B\,$ is bounded uniformly in $\,\mathbb C_+\,$ and this property guarantees (see [6]) relative smoothness of $\,A\,$ with respect to selfadjoint operator $\,L_0,\,$ i.e.
\begin{equation}
\sup_{\|\varphi\|=1,\,\varepsilon\geqslant0}\int_{\mathbb
R}\|AR_0(\tau\pm
i\varepsilon)\varphi\|^2\,d\tau\,<\,\infty\,.
\end{equation}
Further, under the imposed assumptions Jost function $\,e(k)\,$ is bounded away from zero in $\,\mathbb C_+.\,$
Taking this fact into account and applying Lemma 3 once again we come to the conclusion that operator $\,AR_V(k^2)B\,$
is bounded uniformly in $\,\mathbb C_+\,$ and hence (7) and (8) imply (6).

\bigskip

{\bf 5. Wave operators and similarity criterion}

\bigskip

Stationary wave operators to be constructed below intertwine the resolvents $\,R_0(\lambda)\,$ and $\,R_V(\lambda)\,$ and therefore the same holds true for the corresponding propagators $\,e^{itL_0}\,$ and $\,e^{itL_V}.\,$ This enables one to take advantage of time-dependent scattering theory technique in our setting (cf. [12]).

\bigskip

{\bf Theorem 2.} {\it Suppose that bounded potential $\,V\,$ satisfies condition $(\ast)$ and operator $\,L_V=L_0+V\,$ has neither eigenvalues nor spectral singularities. Then direct and inverse wave operators
\begin{eqnarray*}
\Omega_{\pm}&=&s-\lim_{t\to\pm\infty}e^{itL_V}e^{-itL_0},\\
\widetilde{\Omega}_{\pm}&=&s-\lim_{t\to\pm\infty}e^{itL_0}e^{-itL_V}
\end{eqnarray*}
exist and implement the similarity $\,L_V\,=\,\Omega_{\pm}\,L_0\,\widetilde{\Omega}_{\pm}.$ }

\bigskip

{\bf Proof.} Following [6] define operators $\,W_{\pm}\,$ by means of their bilinear forms
\begin{equation}
(W_{\pm}\varphi,\psi)\,=\,(\varphi,\psi)\,\mp\frac1{2\pi
i}\int_{-\infty}^{\infty}\!(AR_0(\tau\pm
i0)\varphi,BR_V(\tau\mp i0)^*\psi)\,d\tau\,.
\end{equation}
For arbitrary $\,\varphi,\psi\in\cal H\,$ vector functions $\,AR_0(\lambda)\varphi\,$ and
$\,BR_V(\lambda)^*\psi\,$ belong (by Statement 2) to Hardy classes $\,{\rm H}_2^+\,$ and ${\rm H}_2^-\,$
respectively. Thus integral on the right-hand-side of (9) represents a bounded linear functional which determines a bounded everywhere defined linear operator.
An intertwining relationship for the resolvents
$$
W_{\pm}\,R_0(\lambda)\,=\,R_V(\lambda)\,W_{\pm}\,
$$
can be verified straightforwardly.
To this end one should match the corresponding forms
$\,(W_{\pm}R_0(\lambda)\varphi,\psi)\,$ and
$\,(W_{\pm}\varphi,R_V(\lambda)^*\psi)\,$ making usage of Hilbert identity and boundary value properties of functions from Hardy classes  (cf. Lemma 2.4 from [6]).

Further introduce operators $\,\widetilde{W}_{\pm}\,$ via the corresponding bilinear forms
\begin{equation*}
(\widetilde{W}_{\pm}\varphi,\psi)\,=\,(\varphi,\psi)\,\pm\frac1{2\pi
i}\int_{-\infty}^{\infty}\!(AR_V(\tau\pm
i0)\varphi,BR_0(\tau\mp i0)^*\psi)\,d\tau\,.
\end{equation*}
In the same way as above Statement 2 guarantees that operators $\,\widetilde{W}_{\pm}\,$ are well-defined and bounded. Moreover it turns out  that $\,\widetilde{W}_{\pm}\,=\,W_{\pm}^{-1};\,$ to verify this fact it suffices to apply the scheme of the proof of Lemma 2.5 from [6]. Thus the relationship
$$
R_V(\lambda)\,=\,{W}_{\pm}\,R_0(\lambda)\,\widetilde{W}_{\pm}\,,
$$
holds true and implies similarity $\,L_V\,=\,W_{\pm}\,L_0\,\widetilde{W}_{\pm}\,$ which certainly involves the inclusion
$\,\widetilde{W}_{\pm}D(L_0)\subset D(L_0)\,.$

Finally let us show that $\,W_{\pm}=s\,$-$\!\!\lim\limits_{t\to\pm\infty}e^{itL_V}e^{-itL_0}\,$ and
\begin{equation}
e^{itL_V}\,W_{\pm}\,\,=\,\,W_{\pm}\,e^{itL_0}\,.
\end{equation}
Intertwining relationship (10) can be derived from the corresponding equality for resolvents under the action of Laplace transform and moreover (10) may be treated as a definition of the propagator $\,e^{itL_V}.\,$ Application of Parseval equality to (9) and taking  (10) into account for arbitrary $\,\varphi,\psi\in\cal H\,$ give
$$
(W_+\varphi,\psi)\,=\,(e^{itL_V}e^{-itL_0}\varphi,\psi)\,\,+\,\,i\int_t^{\infty}(Ae^{-isL_0}\varphi,Be^{-isL_V^*}\psi)\,ds\,.
$$
Consequently the inequality
\begin{equation*}
\big|((W_+-e^{itL_V}e^{-itL_0})\varphi,\psi)\big|\,\leqslant\,\left(\int_t^{\infty}\|Ae^{-isL_0}\varphi\|^2ds\right)^{1/2}
\left(\int_t^{\infty}\|Be^{-isL^*_V}\psi\|^2ds\right)^{1/2}
\end{equation*}
is valid where
$$
\int_0^{\infty}\|Ae^{-isL_0}\varphi\|^2ds\,=\,\frac1{2\pi}\int_{-\infty}^{\infty}
\|AR_0(\tau+i0)\varphi\|^2d\tau\,<\,\infty
$$
by the estimate (8). Besides in virtue of Statement 2 one has
$$
\sup_{\|\psi\|=1}\int_0^{\infty}\!\!\|Be^{-isL^*_V}\psi\|^2ds\,=\,\frac1{2\pi}
\sup_{\|\psi\|=1}\int_{-\infty}^{\infty}\!\!
\|BR_V(\tau-i0)^*\psi\|^2d\tau\,<\,\infty
$$
and therefore $\,\,e^{itL_V}e^{-itL_0}\varphi\to W_+\varphi\,$ as $\,t\to+\infty.\,\,$
Similarly it proves that $\,\,W_-\varphi=\lim\limits_{t\to-\infty}e^{itL_V}e^{-itL_0}\varphi\,\,$ and
$\,\,\widetilde{W}_{\pm}\varphi=\lim\limits_{t\to\pm\infty}e^{itL_0}e^{-itL_V}\varphi\,\,$ for any $\,\varphi\in{\cal H}.$

\bigskip

{\bf Corollary 3.} {\it Under the assumptions of the above theorem propagator $\,e^{itL_V}=W_{\pm}e^{itL_0}\widetilde{W}_{\pm}\,$ is bounded while the corresponding generator $\,L_V\,$ is a spectral operator in the sense of Dunford {\rm [13]}.}

\bigskip

{\bf 6. Acknowledgment}

\medskip

This paper was conceived during the Conference on Semigroups of Operators held in Bedlewo in October, 2013. I would like to express my gratitude to the organizers for their kind invitation and great hospitality.

\bigskip

\begin{center}
\bf REFERENCES
\end{center}

\begin{enumerate}
\item
Naimark M. A. \, Investigation of the spectrum and the expansion in eigenfunctions of a nonselfadjoint differential operator of the second order on a semi-axis // Proc. Mos. Math. Soc., 1954, V.3, P.181-270.
\item
Dolph C. L. \, Recent developments in some nonselfadjoint problems of mathematical physics // Bull. Amer. Math. Soc., 1961, V.67, №1, P.1-69.
\item
Glazman I. M. \, Direct methods of qualitative spectral analysis of singular differential operators, Israel Prog. Scientific Transl., 1965.
\item
Schwartz J. \, Some nonselfadjoint operators // Comm. Pure Appl. Math., 1960, V.13, P.609-639.
\item
Simon B. \, Resonances in one dimension and Fredholm determinants // J. Func. Anal., 2000, V.178, P.396-420.
\item
Kato T. \, Wave operators and similarity for some nonselfadjoint operators // Math. Ann. 1966. V.162. P.258-279.
\item
Stankevich I. V. \, On linear similarity of certain nonselfadjoint operators to selfadjoint operators and on the asymptotic behavior for $t\to\infty$ of the solution of a non-stationary Schr\"odinger equation // Sbornik Math., 1966, V.69, №2, P.161-207.
\item
Stepin S. A. \, The Rayleigh hydrodynamical problem: a theorem on eigenfunction expansion and the stability of plane-parallel flows // Izv. RAN. Ser. Mat., 1996, V.60, №6, P.201-221.
\item
Stepin S. A. \, On the Friedrichs model in one-velocity transport theory // Funct. Anal. Appl., 2001, V.35, №2, P.154-157.
\item Stepin S. A. \,
The Birman-Schwinger principle and Nelkin's conjecture in neutron transport theory // Doklady Math., 2001, V.64, №2, P.152-155.
\item
Stepin S. A. \, Disspative Schr\"odinger operator without a singular continuous spectrum // Sbornik Math., 2004, V.195, №6, P.897-915.
\item
Stepin S. A. \, Wave operators for the linearized Boltzmann equation in one-speed transport theory // Sbornik Math., 2001, V.192, №1, P.141-162.
\item
Dunford N. \, A survey of the theory of spectral operators // Bull. Amer. Math. Soc., 1958, V.64, P.217-274.
\end{enumerate}

\end{document}